\documentclass{article}
\usepackage{fullpage,amssymb,amsmath,amsthm,amsfonts}

\def\eqbd{\mathop{{:}{=}}}

\def\C{\mathbb{C}}
\def\Z{\mathbb{Z}}
\def\F{\mathbb{F}}
\def\H{\mathbb{H}}

\def\ee{{\rm e}}
\def\om{{\rm emax}} 

\def\DD{{\bf\Theta}}
\def\D{{\mathbf\Delta}}
\def\pre{\hbox{\small\sc Pre}}
\def\post{\hbox{\small\sc Post}}

\newtheorem{theorem}{Theorem}[section]

\newtheorem{lemma}[theorem]{Lemma}

\theoremstyle{definition}
\newtheorem{definition}[theorem]{Definition}
\newtheorem{example}[theorem]{Example}

\newtheorem{remark}[theorem]{Remark}

\newcommand{\Sym}{{\operatorname{Sym}}}

\newcommand{\Fr}[1]{{\|#1\|_{\mathrm{F}}}}
\newcommand{\Err}{{\mathrm{Err}}}
\newcommand{\ve}{\varepsilon}

\title{Fast matrix multiplication is stable}
\author{James Demmel\thanks{Mathematics Department and CS Division,
University of California, Berkeley, CA 94720. 
The author acknowledges support of NSF under grants CCF-0444486, ACI-00090127, 
CNS-0325873 and of DOE under grant DE-FC02-01ER25478.},
Ioana Dumitriu\thanks{Mathematics Department, University of California, Berkeley,
CA 94720. The author acknowledges support of the Miller Institute for Basic Research 
in Science.},
Olga Holtz\thanks{Mathematics Department, University of California, Berkeley,
CA 94720.},
and Robert Kleinberg\thanks{Computer Science Division, University of California, 
Berkeley, CA 94720. On leave from the Computer Science Department, Cornell
University. Supported by an NSF Mathematical Sciences Postdoctoral Research
Fellowship. } }
\date{December 6, 2006}
\begin{document}
\maketitle

\begin{abstract}
We perform forward error analysis for a large class of recursive matrix multiplication 
algorithms in the spirit of [D. Bini and G. Lotti, Stability of fast algorithms for 
matrix multiplication, Numer. Math. 36 (1980), 63--72]. 
As a consequence of our analysis, we show that the exponent of matrix multiplication 
(the optimal running time)  can be achieved by numerically stable algorithms. 
We also show that new group-theoretic
algorithms proposed in [H. Cohn, and C. Umans, A group-theoretic approach to fast matrix 
multiplication, FOCS 2003, 438--449] and [H. Cohn, R. Kleinberg, B. Szegedy and 
C. Umans, Group-theoretic algorithms for matrix multiplication, FOCS 2005, 379--388] are all
 included in the class of algorithms to which our analysis applies, and are 
therefore numerically stable.  We perform detailed error 
analysis for three specific fast group-theoretic algorithms. 
\end{abstract}

\section{Introduction}
\label{sec_Intro}

Matrix multiplication is one of the most fundamental operations in numerical linear 
algebra. Its importance is magnified by the number of other problems (e.g., computing 
determinants, solving systems of equations,  matrix inversion, LU decomposition, etc.) 
that are reducible to it (see~\cite[Chapter 16]{BCSbook}).

Starting from Strassen's result~\cite{Strassen} that two square $n{\times}n$ matrices
can be multiplied in $O(n^{2.81})$ operations, a sequence of improvements was made
to achieve ever better bounds on the {\em exponent of matrix multiplication,\/} which 
is the smallest real number $\omega$ for which $n{\times}n$ matrix multiplication can 
be performed in  $O(n^{\omega+\eta})$ operations for any $\eta>0$.  
The complexity of the fastest known method to date due to D.~Coppersmith and 
S.~Winograd~\cite{CoppersmithWinograd} is  about  $O(n^{2.38})$. A new approach based 
on group-theoretic methods was recently developed in~\cite{CU} and~\cite{CKSU}, along
with several ideas that can potentially reduce the bound on $\omega$ to $\omega=2$
(obviously, $\omega$ cannot fall below $2$, since $O(n^2)$ operations are required just
to read off the entries of the resulting matrix).

Along with computational cost,  numerical stability is an equally important factor
for the implementation of any algorithm, since accumulation and propagation of roundoff 
errors may otherwise render the algorithm useless. 
It is the purpose of this
work to analyze recursive fast matrix multiplication
algorithms generalizing Strassen's algorithm, as well as 
the new class of algorithms 
described in~\cite{CU} and~\cite{CKSU}, from
the stability point of view.  
The rounding error analysis of Strassen's method was initiated by Brent
(\cite{brent3,higham90b}, \cite[chap. 23]{higham96}).
In our analysis, we rely on earlier work by Bini and 
Lotti~\cite{BiniLotti}. These results do not apply directly to our setup, 
because they do not account for errors from multiplicative constants,
for nonstationarity in subdividing matrices, and for additional pre- and
post-processing operations (which appear in the methods considered here 
but not in Strassen's method). However, we are able to refine the approach 
of Bini and Lotti to build a  sufficiently inclusive  framework, within which 
the new algorithms proposed in~\cite{CU} and \cite{CKSU} can be analyzed 
in detail.   Combining this framework with a result of Raz \cite{Raz} 
also allows us to prove that there exist
numerically stable matrix multiplication algorithms which
perform $O(n^{\omega+\eta})$ operations, for arbitrarily
small $\eta>0$.

The definition of stability used in this paper measures errors normwise
-- see inequality~(\ref{gen_bound}). This is weaker than the componentwise bound 
satisfied by conventional matrix multiplication \cite[eqn. 3.13]{higham96}.
In fact, Miller~\cite{Miller75} showed that any algorithm
satisfying the componentwise bound must do at least
$n^3$ arithmetic operations. The impact of a normwise
error bound on other algorithms depending on matrix multiplication
has been investigated in~\cite{DemmelHigham92,higham96}.
In~\cite{new} we take up the question whether other linear-algebraic 
algorithms exist that are ``stable'' in some sense and just as fast as the
algorithms considered in this paper. 

This paper is organized as follows: In Section~\ref{sec_aa} 
we discuss the model of arithmetic
and algorithms that is used in the rest of the paper, 
along with some basics on forward error bounds.
In Section~\ref{sec_analysis} we introduce and analyze 
from the stability point of view a wide class 
of recursive algorithms for matrix multiplication.
We begin by discussing Strassen-like algorithms, 
based on recursive partitioning of matrices
into the same number of blocks. 
We prove that all such algorithms are stable, and that this class
of algorithms contains algorithms with running time $O(n^{\omega+\eta})$
for arbitrarily small positive $\eta$.
We then generalize our analysis to algorithms where the number of blocks 
depends on the level of recursion, and finally to algorithms involving 
additional preprocessing before and postprocessing after partitioning 
into blocks.   In Section~\ref{sec_grtheory} we use this
approach to analyze the group-theoretic algorithms from~\cite{CU} and  \cite{CKSU}. 
In particular, in Section~\ref{sec_examples} we perform detailed  stability analysis for
three specific classes of algorithms introduced in~\cite{CU} and  \cite{CKSU}.
The paper ends with a brief discussion of other fast linear algebra algorithms in 
Section~\ref{sec_addla}. 

\section{Model of arithmetic and algorithms}  \label{sec_aa}

We adopt the classical model of rounded arithmetic, where each arithmetic operation 
introduces a small multiplicative error, i.e., the computed value of each arithmetic 
operation $op(a,b)$ is given by $op(a,b)(1+\theta)$ where $|\theta|$ is bounded by
some fixed {\em machine precision\/} $\varepsilon$ but is otherwise arbitrary. The arithmetic
operations in classical arithmetic are $\{+,-, \cdot \}$. 
All of the analysis in this paper applies to matrices with
either real of complex entries, i.e. we interpret the
operands in these arithmetic operations as being either
real or complex numbers.
We assume that the roundoff
errors are introduced by {\em every execution\/} of any arithmetic operation (in contrast
to~\cite{BiniLotti}, where it is assumed that multiplication by entries of the auxiliary coefficient matrices
$U$, $V$ and $W$ is error-free). We further assume that all algorithms output the exact value 
in the absence of roundoff errors (i.e., when all errors $\theta$ are zero). 

For simplicity, let us denote by $\DD$ the set of all errors $\theta$ bounded by 
$\varepsilon$ and by $\D$ the set of all sums $\{1+\theta : \theta \in \DD \}$. We use 
the standard notation
$$ {\bf A}+ {\bf B} \eqbd   \{ a +  b : a\in {\bf A}, \; b\in {\bf B} \}, \quad 
 {\bf A}- {\bf B}  \eqbd   \{ a - b : a\in {\bf A}, \; b\in {\bf B} \},  \quad
 {\bf A}\cdot {\bf B}  \eqbd   \{ a \cdot  b : a\in {\bf A}, \; b\in {\bf B} \} $$ 
for the algebraic sum/difference/product of two sets ${\bf A}$ and ${\bf B}$. 
We will also use the notation ${\bf A}^j$ for the set 
$\underbrace{ {\bf A}\cdot {\bf A} \cdots {\bf A}}_{j\; terms}$. Note that the
error sets $\D^j$ are ordered by inclusion: $\D^j\subseteq \D^{j+1}$ for all $j\in \Z_+$.

We now state the most basic error bound that will be used repeatedly throughout this paper.
Suppose a branch-free algorithm performs a number of arithmetic operations to 
compute a polynomial $f$ in the inputs $x=(x_1, \ldots, x_n)$. Then the resulting 
computed value $f_{comp}$ is a function of both $x$ and the errors $\theta$. Moreover,  
\begin{equation}
f_{comp}(x)\in \sum_j  f_j(x) \D^j \label{sets}
\end{equation} for some polynomials $f_j$ in $x$. 
Suppose $\nu$ is the maximum of the exponents $j$ occurring in the terms $\D^j$ and
$m$ is the number of summands in the expression for $f_{comp}(x)$.  
If the algorithm outputs the correct value $f(x)$ in the absence of roundoff errors, 
then~(\ref{sets}) implies
\begin{equation}
 |f_{comp}(x)-f(x)| \leq m \max_j |f_j(x)|( (1+\varepsilon)^\nu -1)=m\max_j |f_j(x)| \nu
\varepsilon +O(\varepsilon^2). \label{poly_error}
\end{equation}

In particular, suppose that addition of $n$ quantities $x_1, \ldots, x_n$ is performed by 
running the classical parallel $\lceil \log_2 n \rceil$-step algorithm as a straight-line algorithm,
i.e., computing the total sum by adding the two sums $\sum_{j=1}^{\lceil n/2  \rceil} 
x_j$ and $\sum_{j=\lceil n/2 \rceil +1}^n x_j$, and recursively computing each of the two
sums in the same manner. Then the resulting computed value $\Sigma_{comp}(x)$ 
lies in the set $$ \sum_{j=1}^n  x_j \D^{\lceil  \log_2 n \rceil},$$ hence
\begin{equation}| \Sigma_{comp} (x)- \sum_j x_j| \leq \max_j |x_j| \lceil \log_2 n \rceil
\varepsilon +O(\varepsilon^2).  \label{sum_error}  
\end{equation}

\section{Error analysis for recursive matrix multiplication algorithms} \label{sec_analysis}

In this section, we perform forward error analysis for three 
classes of recursive matrix multiplication algorithms, 
starting with the Strassen-like algorithms based on stationary
partitioning, then generalizing to algorithms with non-stationary partitioning,  
and finally to the algorithms of the kind developed in~\cite{CU} and~\cite{CKSU}. 

The error analysis in Sections~\ref{sec_stat} and~\ref{sec_nonstat} 
is done with respect to the entry-wise maximum norm on $A$, $B$, $C=AB$, 
while the analysis in Section~\ref{sec_extra}
is for an arbitrary matrix norm satisfying an extra monotonicity assumption. 
All our bounds are of the form 
\begin{equation} 
\| C_{comp}-C\|\leq \mu(n) \varepsilon \|A\| \, \|B \| +O(\varepsilon^2), \label{gen_bound} \end{equation}
with $\mu(n)$ typically low degree polynomials in the order $n$ of the 
matrices involved, so that $\mu(n)=O(n^c)$ for some constant $c$. 
Note that one can easily switch from one norm 
to another at the expense of picking up  additional factors that will 
depend $n$, using the equivalence of norms on a finite-dimensional space.

Later, in Section~\ref{sec_examples}, we will give values of the exponent $c$
for sample algorithms. Here, we argue that the exact value of $c$ 
does not greatly impact the complexity of practical matrix multiplication, 
in the sense of the bit-complexity for computing
$AB$ to a desired accuracy 
$\| C_{comp} - C \| \leq \epsilon_0 \|A\| \cdot \|B\|$,
for a given $\epsilon_0 < 1$.

Since $\|C\|$ can be about as large as $\|A\| \cdot \|B\|$, 
the bound~(\ref{gen_bound}) is interesting only 
when $\mu(n)\varepsilon < 1$.   Any algorithm will have $c\geq 1$, 
since even just straightforwardly computing a dot product has $c=1$. 
Thus $n^c \varepsilon\leq 1$, so $n\varepsilon \leq 1$,
and $b\eqbd \log_2 (1/\varepsilon)\geq \log_2 n$, 
where $b$ is the number of bits used to 
represent the fractional part of the floating point numbers. 

Now suppose we want to choose $\epsilon$ just small enough ($b$ just
large enough) to guarantee $\| C_{comp} - C \| \leq \epsilon_0 \|A\| \cdot \|B\|$,
and ask how the complexity of the resulting algorithm depends on the
exponent $c$. Setting $\epsilon_0 = n^c \epsilon$, we get
$b = \log_2 (1/\epsilon) = \log_2 (1/ \epsilon_0) + c \cdot \log_2 n \geq \log_2 n$,
i.e. the number of bits $b$ needed grows proportionally to $\log_2 n$.
The cost of $b$-bit arithmetic is in the range from
$O(b^2)$ (done straightforwardly) down to $O(b^{1+o(1)})$ 
(using Sch\"onhage-Strassen~\cite{SchonhageStrassen}). 
Therefore, the bit-complexity of computing the answer with error 
proportional to $\epsilon_0$ will be at most 
a polylog($n$) factor larger than the bound $O(n^{\omega})$ 
gotten from ignoring bit-complexity, and only slightly superlinearly 
(up to quadratically) dependent on $c$.

\subsection{Stationary partition algorithms}  \label{sec_stat}

We next recall some basic notions related to recursive matrix multiplication 
algorithms.  This section is closely related to the paper~\cite{BiniLotti} by Bini and Lotti. 
However, our approach is more inclusive, as will be explained later in this section as we 
develop pertinent details.

We consider recursive algorithms for matrix multiplication. A {\em bilinear non-commutative
algorithm\/}~(see \cite{BiniLotti} or \cite{BrockettDobkin}) that computes products of $k\times k$
 matrices $C=AB$ over a ground field $\F$ using $t$ non-scalar multiplications is determined by 
three  $k^2{\times}t$ matrices $U$, $V$ and $W$ with elements in a subfield $\H \subseteq \F$ 
such that
\begin{equation} 
c_{hl}=\sum_{s=1}^t w_{rs} P_s, \;\; {\rm where} \;\; 
P_s\eqbd \left(\sum_{i=1}^{k^2} u_{is} x_i \right) \left(\sum_{j=1}^{k^2} v_{js} y_j \right),
 \;\;\;\;  r=k(h-1)+l,\; \;\; h,l=1,\ldots,k,   \label{main} 
\end{equation}
where $x_i$ (resp. $y_j$) are the elements of $A=(a_{ij})$ (resp. of $B=(b_{ij})$) ordered column-wise, 
and $C=(c_{ij})$ is the product $C=AB$.  

For an arbitrary $n$, the algorithm consists in recursive partitioning and using
formula~(\ref{main}) to compute products of resulting block matrices. More precisely,
suppose that $A$ and $B$ are of size $n{\times}n$, where $n$ is a power of $k$ (which
can always be achieved by augmenting the matrices $A$ and $B$ by zero columns and rows).
Partition $A$ and $B$ into $k^2$ square blocks $A_{ij}$, $B_{ij}$ of size $(n/k){\times}(n/k)$.
Then the blocks $C_{hl}$ of the product $C=AB$ can be computed by applying~(\ref{main}) to 
the blocks of $A$ and $B$, where each block $A_{ij}$ , $B_{ij}$ has to be again partitioned 
into $k^2$ square sub-blocks to compute the $t$ products $P_s$ and then the blocks $C_{hl}$. 
The algorithm obtained by running this recursive procedure $\log_k n$ times computes the
product $C=AB$ using at most $O(n^{\log_k t})$ multiplications.  

Now we are in a position to analyze recursive matrix multiplication algorithms.
We first look at the outermost recursion, denoting the blocks of $A$ ordered
column-wise by $X_i$ and the blocks of $B$ ordered column-wise by $Y_j$. We will index 
the levels of recursion by $j=1, \ldots, p\eqbd \log_k n$, increasing as we go down.
Since multiplication by an element of $U$ or $V$ introduces a multiple of $1+\theta$ 
for some $\theta\in\DD$ and  since~(\ref{main}) and~(\ref{sum_error}) hold,  
the computed value $P^{[1]}_{s,comp}$ for each quantity $P^{[1]}_s$ is obtained by running
the fast matrix multiplication algorithm on the obtained pairs of $(n/k){\times}(n/k)$ 
matrices
\begin{equation}
A^{[1]}_{s,comp}\in \sum_{i=1}^{k^2} u_{is} X_i \D^{1+\alpha_s}, \qquad 
B^{[1]}_{s,comp}\in \sum_{j=1}^{k^2} v_{js} Y_j \D^{1+\beta_s}, \label{int1} 
\end{equation}
where $ \alpha_s\eqbd \lceil \log_2 a_s \rceil$, $ \beta_s\eqbd \lceil \log_2 b_s \rceil$, and
$$ a_s\eqbd \# \{ u_{is}\;  : \; u_{is}\neq 0, \; i=1, \ldots, k^2 \}, \quad
b_s\eqbd \# \{ v_{js}\;  : \; v_{js}\neq 0, \; j=1, \ldots, k^2 \} \quad {\rm for} \;\;
s=1, \ldots, t.$$ 
The matrices $A^{[1]}_{s,comp}$, $B^{[1]}_{s,comp}$ are further partitioned and the same procedure
is applied to the obtained blocks, etc., $\log_k n$ times, until the resulting blocks 
all have size $1{\times}1$. To see how the errors propagate, note that if the inputs 
to~(\ref{main}) are given as sums of certain matrices $A_{\phi},B_{\psi},$
each with a possible error in $\D^\alpha$,
$\D^\beta$, respectively (i.e., as elements of the sets $\sum_{A_\phi} A_\phi \D^\alpha$, $\sum_{B_\psi} 
B_\psi \D^\beta$), then the resulting inputs at the next level are elements of the sets
$$ \sum_{A_\phi} \sum_{i=1}^{k^2} u_{is} X_i(A_\phi) \D^{\alpha+\alpha_s+1}, \qquad  
   \sum_{B_\psi} \sum_{j=1}^{k^2} u_{js} Y_j(B_\psi) \D^{\beta+\beta_s+1},\qquad  {\rm respectively.} $$
Thus by going all the way down to the $\log_k n$ level we multiply each element of the original
input matrix $A$ by error terms in $\D^{(1+\alpha_s)\log_k n }$ and by $\log_k n$ elements of $U$.
Likewise, each element of $B$ is multiplied by error terms from the set 
 $\D^{(1+\beta_s)\log_k n }$ and $\log_k n$ elements of $V$. 

Now, at the lowest level of our recursive scheme, we begin to put together quantities
$P^{[p]}_{s,comp}$. The lowest-level computation of $P^{[p]}_{s,comp}$ 
is simply a scalar multiplication, so it brings in an additional factor 
from $\D$. Then the quantities $c_{h_p,l_p,comp}$ are computed by~(\ref{main}).
Note that, in general, 
\begin{equation}
 C_{hl,comp}\in \sum_{s=1}^t w_{rs} P_{s,comp} \D^{\gamma_r+1}, \label{int2}
\end{equation}
where $\gamma_r\eqbd \lceil \log_2 c_r   \rceil$ and
$$ c_r\eqbd \# \{ w_{rs}\;  : \; w_{rs}\neq 0, \; s=1, \ldots, t \}\quad
{\rm for} \;\; r=1, \ldots, k^2.$$
Also, $a_s b_s$ terms of the kind $w_{rs}u_{is}v_{js}X_iY_j$ need to be added to produce 
$P_{s,comp}$ from the products $(X_iY_j)$ using formula~(\ref{int1}). So 
$\sum_{s=1}^t a_s b_s \xi_{rs}$ terms containing a product $x_i y_j$, where $x_i$
is an entry of $A$, $y_j$ is an entry of $B$, are  added to produce $C_{h_pl_p,comp}$ 
by formula~(\ref{int2}), where 
$$ \xi_{rs}\eqbd \left\{ \begin{array}{ll} 1 & w_{rs}\neq 0  \\
0 & w_{rs}=0.  \end{array}  \right. $$
Following~\cite{BiniLotti}, we denote by $\ee$ the vector with components 
$ \ee_r\eqbd  \sum_{s=1}^t a_s b_s \xi_{rs}$, and by $\om$ the maximum
$\om\eqbd \max_r \ee_r$.
 
As the second part of the recursive procedure is run from the bottom to the 
top, we ``assemble'' all the blocks $C_{h_j,l_j,comp}$ from the blocks at
the previous level. When the algorithm terminates, each resulting 
element of $C$ is then determined by the choice 
of block indices $(h_1, l_1)$, $\ldots$, $(h_p,l_p)$  and is the sum of $e_{r_1}
\cdots e_{r_p}$ terms  $c_{h_1,l_1, \ldots, h_p,l_p,comp}$, where $r_q\eqbd (h_q-1)k+l_q$. 
Each term $c_{h_1,l_1, \ldots, h_p,l_p,comp}$ is a product of an element of $A$ and an
element of $B$, an element of $\D^\mu$ where $\mu$ is at most  $1+\max_{r,s}(\alpha_s+\beta_s+
\gamma_r+3)\log_k n$, and  $\log_k n$ elements  of $U$, of $V$ and of $W$. Using the 
maximum-entry norm $\|M\|  \eqbd \max_{ij} |m_{ij}|$ for a matrix $M$,   
we therefore arrive at the bound
\begin{equation}
 \|C_{comp}-C\| \leq (1+\max_{r,s}(\alpha_s+\beta_s+\gamma_r+3)\log_k n) \cdot
( \om \cdot \|U\|  \,\|V\| \, \|W\| )^{\log_k n} \|A\|  \, \|B\| 
 \varepsilon +O(\varepsilon^2).
\label{bound_stat}
\end{equation}

We can now summarize this formally as a theorem. 

\begin{theorem} \label{thm:stationary}
A bilinear non-commutative algorithm for matrix multiplication based
on stationary partitioning is stable. It satisfies the error bound~(\ref{gen_bound}) 
where $\|\cdot\|$ is the maximum-entry norm and where
$$ \mu(n) = (1+\max_{r,s}(\alpha_s+\beta_s+\gamma_r+3)
\log_k n) \cdot ( \om \cdot \|U\|  \,\|V\| \, \|W\| )^{\log_k n}. $$
\end{theorem}
\begin{remark}
Note that in Theorem~\ref{thm:stationary},
$\mu(n) = O \left(
n^{\log_k(\om \cdot \|U\|  \,\|V\| \, \|W\| ) + o(1)} \right).$
This confirms our statement, following equation (\ref{gen_bound}),
that the term $\mu(n)$ in our error bound is polynomial in $n$.
\end{remark}

\begin{theorem} \label{thm:omega-is-stable}
For every $\eta>0$ there exists an algorithm
for multiplying $n$-by-$n$ matrices which performs
$O(n^{\omega+\eta})$ operations (where
$\omega$ is the exponent of matrix multiplication)
and which is numerically stable, in the sense that it 
satisfies the error bound~(\ref{gen_bound}) with 
$\mu(n)=O(n^c)$ for some constant $c$ depending on
$\eta$ but not $n$.
\end{theorem}
\begin{proof}
It is known that the exponent
of matrix multiplication is achieved by bilinear
non-commutative algorithms~\cite{Raz}.  More 
precisely, using the terminology of~\cite{Raz}, 
for any arithmetic circuit of size $S$ which computes 
the product of two input matrices $A$, $B$ over 
a field of characteristic zero, there is 
another arithmetic circuit of size $O(S)$
which also computes the product of $A$ and $B$ 
and is a \emph{bilinear circuit}, meaning that
it has the following structure.  There are
two subcircuits $\mathcal{C}_1$, $\mathcal{C}_2$, 
whose outputs are 
linear functions of the entries $a_{ij}$
(resp. $b_{ij}$).  Then there is one layer
of product gates, each of which multiplies
one output of $\mathcal{C}_1$
with one output of $\mathcal{C}_2$.
Then there is a subcircuit $\mathcal{C}_3$ 
whose inputs are the outputs of these 
product gates, and whose outputs are the
entries of the matrix product.
The only operations performed inside
subcircuits $\mathcal{C}_1$, $\mathcal{C}_2,
\mathcal{C}_3$ are addition and scalar 
multiplication.  Every bilinear circuit 
corresponds to a bilinear noncommutative 
algorithm as expressed in~(\ref{main})
the outputs of $\mathcal{C}_1$, $\mathcal{C}_2$ are
the linear forms $\{\sum u_{is} x_i\}$,
$\{ \sum v_{js} y_j \},$ respectively.
The product gates compute the numbers
$P_s$.  The circuit $\mathcal{C}_3$
computes the linear forms $\sum w_{rs} P_s$.

By the definition of $\omega$, for some constant $C$ 
there exist arithmetic circuits of size less than
$C k^{\omega+\eta/2}$ which compute a $k \times k$
matrix product, for every $k$.  By the preceding 
paragraph, we can assume that these circuits are 
bilinear circuits.  This means that for every
$k$ there exists a bilinear noncommutative algorithm
for $k \times k$ matrix multiplication 
using $t < C k^{\omega + \eta/2}$
non-scalar multiplications.
Choose $k_0$ large enough that
$C k_0^{\omega+\eta/2} < k_0^{\omega+\eta}.$
Using the bilinear non-commutative algorithm
for this value of $k_0$ and applying 
Theorem~\ref{thm:stationary}, we obtain 
the theorem.
\end{proof}

\subsection{Non-stationary partition algorithms}  \label{sec_nonstat}

The analysis from the preceding section generalizes easily to bilinear matrix multiplication 
algorithms based on non-stationary partitioning. In that case, the matrices $A^{[j]}_{s,comp}$ 
and $B^{[j]}_{s,comp}$ are partitioned into $k{\times}k$ square blocks, but $k$ depends on
the level  of recursion, i.e., $k=k(j)$, and the corresponding matrices $U$, $V$ and $W$ also
depend on $j$: $U=U(j)$, $V=V(j)$, $W=W(j)$.   Otherwise the algorithm proceeds exactly like
in the previous section. Suppose such an algorithm applied to $n{\times}n$ matrices requires
$p$ levels of recursion, so that $\prod_{j=1}^p k(j)=n$. For each level $j$, we can define
quantities $\alpha_s(j)$ analogously to $\alpha_s$,  $\beta_s(j)$ analogously to $\beta$,
$\om(j)$ analogously to $\om$. We then use the same reasoning as above to obtain the following 
error bound for  non-stationary partition algorithms:
\begin{eqnarray}
\|C_{comp}-C\| & \leq & (1+\sum_j\max_{r,s}(\alpha_s(j)+\beta_s(j)+\gamma_r(j)+3)) \nonumber \\
 && {\times} \left(  \prod_j\om(j) \|U(j)\| \,\|V(j)\|\, \|W(j)\| \right) \|A\| \, \|B\|
 \varepsilon +O(\varepsilon^2). 
\label{bound_nonstat}
\end{eqnarray}

\begin{theorem} A bilinear non-commutative algorithm for matrix multiplication based
on non-stationary partitioning is stable. It satisfies the error bound~(\ref{gen_bound}) 
where $\| \cdot\|$ is the maximum-entry norm and where
$$ \mu(n) =(1+\sum_j\max_{r,s}(\alpha_s(j)+\beta_s(j)+\gamma_r(j)+3)) 
\left( \prod_j\om(j) \|U(j)\| \,\|V(j)\|\, \|W(j)\| \right). $$
\end{theorem}

\subsection{Algorithms that combine partitioning with pre- and post-processing} \label{sec_extra}

Finally consider algorithms that combine recursive non-stationary partitioning with pre- and
post-processing given by linear maps $\pre_n()$ and $\post_n()$ acting on matrices of an 
arbitrary order $n$.
More specifically, the matrices $A$ and $B$ are each pre-processed, then partitioned into 
blocks, respective pairs of blocks are multiplied recursively and assembled into a large matrix, 
which is then post-processed to obtain the resulting matrix $C$ (see 
Section~\ref{sec_abelian_stp_alg} for concrete examples of pre- and processing operators).

We assume again that the partitioning
is non-stationary, i.e., that at  level $j$ of the recursion all matrices  are of order
$n_j\eqbd \prod_{l\geq j} k(l)$  and are partitioned into $t_j$ blocks of size $n_{j+1}$.
(At the lowest $p$th level of recursion, $n_p=1$, while at the top level of recursion, $n_1=n$, 
i.e., $\prod_{j=1}^p k(j)=n$.)

For this analysis, we will be working with a {\em consistent\/} norm $\| \cdot \|$ 
defined for  matrices of all sizes and satisfying the condition
\begin{equation}
\max_s  \|M_{s} \|\leq  \|M \|\leq \sum_{s}\|M_{s} \|
\label{normcond}
\end{equation}
 whenever the matrix $M$ is partitioned into
 blocks $(M_{s})_s$ (an example of such a norm is provided by $\|\cdot \|_2$). 
Note that the previously used maximum-entry norm satisfies~(\ref{normcond}) but is not 
consistent, i.e., fails to satisfy
$$ \| A B\| \leq \|A \| \cdot \|B \| \qquad \hbox{\rm for all } A, \; B.$$ 

We denote the norms of pre- and post- processing maps subordinate to
the norm $\| \cdot \|$ by $\| \cdot \|_{op}$.
Suppose that the pre- and post-processing is performed with errors 
$$ \| \pre_n (M)_{comp}-\pre_n (M) \| \leq f_{pre}(n) \ve \|M\|+ O(\varepsilon^2), \;\;
\| \post_n (M)_{comp}-\post_n (M) \| \leq f_{post}(n) \ve \|M\|+ O(\varepsilon^2), $$
where $n$ is the order of the matrix $M$. 
As before, we denote by $\mu(n)$ the coefficient of $\varepsilon$ in the final error 
bound~(\ref{gen_bound}). 

The function $\mu$ can be found recursively as follows. Consider one level of recurrence 
where matrices of order $n$ are partitioned into, say, $t$ matrices of order $n/k$. Denote the matrix $\pre_n(A)$ by $\hat{A}$ and $\pre_n(B)$ by $\hat{B}$. The computed matrix
$\hat{A}_{comp}$ ($\hat{B}_{comp}$, resp.) is within $f_{pre}(n) \ve \|A\|$ ($f_{pre}(n) 
\ve \|B\|$, resp.) from $\hat{A}$ ($\hat{B}$, resp.). The matrices $\hat{A}$ and $\hat{B}$ 
are further partitioned, which does not introduce additional errors. Thus
$$ \| \hat{A}_{s,comp} -\hat{A}_{s} \|\leq f_{pre}(n) \ve \|A\| +O(\ve^2),  \quad
 \| \hat{B}_{s,comp} -\hat{B}_{s} \|\leq f_{pre}(n) \ve \|B\| +O(\ve^2).$$
The blocks $\hat{A}_{s,comp}$ and $\hat{B}_{s,comp}$ are then multiplied recursively,
which, for each pair of blocks, introduces an error of size $\mu(n/k)\ve \|\hat{A}_{s,comp} \|\,
 \|\hat{B}_{s,comp} \|$. Denoting the computed products by $\hat{C}_{s,comp}$, we thus obtain
\begin{eqnarray*}
 \| \hat{C}_{s,comp} - \hat{A}_{s,comp} \hat{B}_{s,comp}\|  
& \leq & \mu(n/k) \ve  \|\hat{A}_{s,comp} \| \, \|\hat{B}_{s,comp} \|+O(\ve^2) \\ 
& \leq & \mu(n/k) \ve  \| \hat{A}_{s} \| \, \| \hat{B}_{s} \|+O(\ve^2)  \\
& \leq & \mu(n/k) \ve \| \hat{A} \| \, \| \hat{B} \|+O(\ve^2) \\
& \leq & \mu(n/k) \ve \|\pre_n \|^2_{op}\, \|\ A \| \, \|\ B \|+O(\ve^2).       
\end{eqnarray*}
We now apply the triangle inequality to evaluate $\|\hat{C}_{s,comp}-\hat{A}_{s} \hat{B}_{s}\|$.
Rewriting $\hat{A}_{s}$ as $\hat{A}_{s,comp}+(\hat{A}_{s}- \hat{A}_{s,comp})$ and $\hat{B}_{s}$
as $\hat{B}_{s,comp}+(\hat{B}_{s}- \hat{B}_{s,comp})$, we get
\begin{eqnarray*} 
\|\hat{C}_{s,comp}-\hat{C}_{s} \|& \leq &\| \hat{C}_{s,comp}-\hat{A}_{s,comp} 
\hat{B}_{s,comp} \|+ \|\hat{A}_{s,comp} \hat{B}_{s,comp} -\hat{A}_{s} \hat{B}_{s} \| \\  
&\leq & \mu(n/k) \ve  \| \pre_n\|^2_{op} \, \|A\| \, \|B\| + 2 f_{pre}(n) \ve \|A \|\, \| B\| +O(\ve^2). 
\end{eqnarray*}
Summing up over all $s=1, \ldots, t$ and taking into account the assumed properties of the norm 
$\| \cdot\|$, we therefore obtain
$$ \| \hat{C}_{comp} -\hat{C} \| \leq \sum_s \| \hat{C}_{s,comp}- \hat{C}_s \|
\leq t \left(\mu(n/k) \ve  \| \pre_n\|^2_{op} \, \|A\| \, \|B\| + 2 f_{pre}(n) \ve \|A \|\, \| B\|\right) +O(\ve^2).$$ 
Finally, the post-processing step rescales the obtained error by the norm $\|\post_n \|_{op}$ and 
adds another error term of order 
$$ f_{post}(n) \| \hat{C}_{comp} \| \ve \leq  f_{post}(n) \|\pre_n \|_{op}^2\, \| A \|\,  \|B\| +O(\ve^2).$$ 
Altogether, this gives the recurrence
$$  \mu(n) = \mu(n/k) t \|\post_n \|_{op} \, \|\pre_n \|_{op}^2 + 2 f_{pre}(n) t \|\post_n\|_{op} +f_{post}(n) \|\pre_n \|_{op}^2.  $$ 
The same argument is applicable to each level $j$ of recurrence, with $n$ replaced by $n_j$,
$n/k$ replaced by $n_{j+1}$, and $t$ replaced by $t_j$. We now state this formally as a theorem.

\begin{theorem} \label{thm_crude}
Under the assumptions of this section, a recursive matrix multiplication algorithm based on 
non-stationary  partitioning with pre- and post-processing is stable. It satisfies the error 
bound~(\ref{gen_bound}), with the function $\mu$ satisfying the recursion
$$  \mu(n_j) = \mu(n_{j+1}) t_j \|\post_{n_j} \|_{op} \, \|\pre_{n_j} \|_{op}^2 + 2  f_{pre}(n_j) t_j \|\post_{n_j}\|_{op} 
+f_{post}(n_j) \|\pre_{n_j} \|_{op}^2  \,  , \qquad j=1, \ldots, p.  $$

\end{theorem}

\section{Group-theoretic recursive algorithms} \label{sec_grtheory}

To perform the error analysis of the group-theoretic matrix multiplication 
algorithms defined in~\cite{CU} and~\cite{CKSU}, we must first recall some definitions
and facts about those algorithms.  The relevant material is
reviewed in Section~\ref{sec_grtheory_background}.
In Section~\ref{sec_abelian_stp_alg}
we define the class of group-theoretic algorithms ---
called \emph{abelian simultaneous triple product (abelian
STP) algorithms} --- and we introduce a 
running example (i.e. a specific algorithm in this
class) for the purpose of concreteness.
This class of algorithms encompasses all of the fast matrix 
multiplication algorithms described in~\cite{CKSU}, and is a special 
case of the group-theoretic algorithms defined in~\cite{CU}.  
We refer the reader to~\cite{CKSU-companion} for a proof
of the correctness of these algorithms, as well as an
analysis of their running time.
In Section~\ref{sec_grtheory_analysis} we will
apply the analysis from Section~\ref{sec_analysis}
to derive error bounds for abelian STP algorithms.
In Section~\ref{sec_examples} we will cite
some specific examples of such algorithms and
evaluate their error bounds.

\subsection{Background material} \label{sec_grtheory_background}

We begin by recalling some basic definitions
from algebra.
\begin{definition}[semidirect product]
If $H$ is any group and $Q$ is a group which
acts (on the left) by automorphisms of $H$, with 
$q \cdot h$ denoting
the action of $q \in Q$ on $h \in H$, then 
the \emph{semidirect product} $H \rtimes Q$ is the
set of ordered pairs $(h,q)$ with the
multiplication law
\begin{equation} \label{eqn:semidirect}
(h_1,q_1) (h_2,q_2) = (h_1(q_1 \cdot h_2), q_1 q_2).
\end{equation}
We will identify $H \times \{1_Q\}$ with $H$
and $\{1_H\} \times Q$ with $Q$, so that 
an element $(h,q) \in H \rtimes Q$ may
also be denoted simply by $hq$.  Note that
the multiplication law of $H \rtimes Q$
implies the relation $qh = (q \cdot h) q$.
\end{definition}

\begin{definition}[wreath product]
If $H$ is any group, $S$ is any finite set,
and $Q$ is a group with a left action on $S$,
the \emph{wreath product} $H \wr Q$ is the semidirect
product $(H^S) \rtimes Q$ where $Q$ acts on 
the direct product of $|S|$ copies of $H$
by permuting the coordinates according to 
the action of $Q$ on $S$.  (To be more
precise about of the action of $Q$ on
$H^S$, if an element $h \in H^S$ is represented
as a function $h: S \rightarrow H$, then 
$q \cdot h$ represents the function 
$s \mapsto h(q^{-1}(s)).$)
\end{definition}

\begin{example}[running example, part 1]
\label{ex_part1}
Throughout this section, we will work with a 
running example of an abelian STP algorithm
based on a specific finite abelian group $H$ with
$4096$ elements, and its wreath product with a
two-element group.  Consider the
set $S \eqbd \{0,1\}$ and a two-element group $Q$ whose
non-identity element acts on $S$
by swapping $0$ and $1$.  Let $H$ be the group
$(\Z/16)^3$ whose elements are ordered triples
of integers $(x_0,x_1,x_2)$ modulo $16$.  An
element of $H^S$ is an ordered pair of elements
of $H$, which can be represented as a $2$-by-$3$
matrix 
$$ \left( \begin{array}{lll}
x_{00} & x_{01} & x_{02} \\ x_{10} & x_{11} & x_{12}
\end{array} \right)$$
of integers modulo $16$.  An element of $H \wr Q$
is an ordered pair $(X,q)$ where $X$ is a matrix
as above, and $q = \pm 1$.  An example of the 
multiplication operation in $H \wr Q$ is given
by the formula for $(X,-1) \cdot (Y,-1)$:
$$ \left( \left( \begin{array}{lll}
x_{00} & x_{01} & x_{02} \\ x_{10} & x_{11} & x_{12}
\end{array} \right), -1 \right)
\cdot
\left( \left( \begin{array}{lll}
y_{00} & y_{01} & y_{02} \\ y_{10} & y_{11} & y_{12}
\end{array} \right), -1 \right) = 
\left( \left( \begin{array}{lll}
x_{00} + y_{10} & x_{01} + y_{11} & x_{02} + y_{12} \\ 
x_{10} + y_{00} & x_{11} + y_{01} & x_{12} + y_{02}
\end{array} \right), 1 \right).
$$
Notice that the rows of $Y$ were swapped before adding
it to $X$.

An alternative description of $H \wr Q$ is that it has
generators $a_0,a_1,b_0,b_1,c_0,c_1$ satisfying the 
following relations:
\begin{enumerate}
\item  $a_0,a_1,b_0,b_1,c_0,c_1$ collectively generate
the group $(\Z/16\Z)^6$.
\item  $q^2$ is the identity element.
\item  $q a_0 = a_1 q, \, q b_0 = b_1 q, \, q c_0 = c_1 q.$
\end{enumerate}
\end{example}

\begin{example}
This example generalizes the preceding one.
When $S$ is the set $\{1,2,\ldots,n\}$,
we use the notation $\Sym_n$ to denote the 
group of all permutations of $S$, acting on $S$ 
in the obvious way, i.e. $\pi \cdot s = \pi(s).$
Each element of the wreath product $H \wr \Sym_n$ 
may be uniquely represented as a product
$h \pi$ where $\pi \in \Sym_n$ and 
$h = (h_1,h_2,\ldots,h_n) \in H^n$.  
The multiplication law of $H \wr \Sym_n$
is given by the formula:
\begin{equation} \label{eqn:wreath}
(h \pi)(h' \pi') = 
h (\pi \cdot h') \pi \pi' =
\left( h_1 h'_{\pi^{-1}(1)},
h_2 h'_{\pi^{-1}(2)},
\ldots,
h_n h'_{\pi^{-1}(n)} \right)
\pi \pi'.
\end{equation}
\end{example}
\bigskip
Next we recall some definitions and theorems
from~\cite{CU} and~\cite{CKSU}.
If $S,T$ are subsets of a group $G$, we use the notation
$Q(S,T)$ to denote their right quotient set, i.e.
$$Q(S,T) \eqbd \{s t^{-1} \,:\, s \in S, t \in T\}.$$
We use the notation $Q(S)$ as shorthand for $Q(S,S)$.
\begin{definition}[triple product property,
simultaneous triple product property]
If $H$ is a group and $X,Y,Z$ are three subsets,
we say $X,Y,Z$ satisfy the \emph{triple product
property} if it is the case that for all
$q_x \in Q(X), q_y \in Q(Y), q_z \in Q(Z),$ if
$ q_x q_y q_z = 1 $
then $q_x=q_y=q_z=1.$  

If $\{(X_i,Y_i,Z_i) \,:\, i \in I\}$ is a 
collection of ordered triples of subsets of 
$H$, we say that this collection satisfies
the \emph{simultaneous triple product property (STPP)}
if it is the case that for all $i,j,k \in I$ and all
$q_x \in Q(X_i,X_j),
q_y \in Q(Y_j,Y_k),
q_z \in Q(Z_k,Z_i),$ if
$q_x q_y q_z = 1$ then
$q_x = q_y = q_z = 1$ and $i=j=k$.
\end{definition}

\begin{example}[running example, part 2]
\label{ex_part2}
In our running example, the group $H$ is $(\Z/16\Z)^3$.
Consider the following three subgroups of $H$.
\begin{eqnarray*}
X & \eqbd & (\Z/16\Z) \times \{0\} \times \{0\} \\
Y & \eqbd & \{0\} \times (\Z/16\Z) \times \{0\} \\
Z & \eqbd & \{0\} \times \{0\} \times (\Z/16\Z) \\
\end{eqnarray*}
We claim that $X,Y,Z$ satisfy the triple product property.
Since $H$ is an abelian group, we will denote the group
operation and the identity element  using additive notation
rather than multiplicative notation.  Thus the triple
product property is the assertion that if $q_x \in Q(X),
q_y \in Q(Y), q_z \in Q(Z),$ and $q_x + q_y + q_z = 0$,
then $q_x = q_y = q_z = 0.$
Note first that $Q(X)=X, Q(Y)=Y, Q(Z)=Z$ because $X,Y,Z$ 
are subgroups.  The elements $q_y, q_z$ have $0$ in their
first component, so the first component of $q_x + q_y + q_z$ is
equal to the first component of $q_x$.  This shows that the
first component of $q_x$ is $0$, which implies that $q_x=0$.
By similar arguments, $q_y=0$ and $q_z=0$, which confirms
the triple product property.

Now consider the following six subsets of $H$:
\begin{eqnarray*}
\overline{X}_0  \eqbd  \{1,2,\ldots,15\} \times \{0\} \times \{0\}
& \qquad &
\overline{X}_1  \eqbd  \{0\} \times \{1,2,\ldots,15\} \times \{0\} \\
\overline{Y}_0  \eqbd  \{0\} \times \{1,2,\ldots,15\} \times \{0\} 
& \qquad &
\overline{Y}_1  \eqbd  \{0\} \times \{0\} \times \{1,2,\ldots,15\} \\
\overline{Z}_0  \eqbd  \{0\} \times \{0\} \times \{1,2,\ldots,15\} 
& \qquad &
\overline{Z}_1  \eqbd  \{1,2,\ldots,15\} \times \{0\} \times \{0\}
\end{eqnarray*}
We claim that $(\overline{X}_0,\overline{Y}_0,\overline{Z}_0)$
and $(\overline{X}_1,\overline{Y}_1,\overline{Z}_1)$ satisfy
the simultaneous triple product property.  Suppose that 
$i,j,k \in \{0,1\}$ and $q_x \in Q(X_i,X_j), q_y \in Q(Y_j,Y_k),
q_z \in Q(Z_k,Z_i).$  Suppose moreover that $q_x + q_y + q_z = 0$.
If $i=j=k$, then we may argue as before that 
$\overline{X}_i,\overline{Y}_i,\overline{Z}_i$ satisfy
the triple product property and therefore $q_x=q_y=q_z=0$
as desired.  If $i,j,k$ are not all equal, we may perform
a case analysis for each of the six possible ordered triples 
$(i,j,k)$, in each case obtaining a conclusion which
contradicts the assumption that $q_x + q_y + q_z = 0$.
We illustrate this by considering the case $i=j=0, k=1.$
In this case $q_x \in Q(X_0,X_0)$ has zero in its second
component, as does $q_z \in Q(Z_1,Z_0).$  But $q_y \in
Q(Y_0,Y_1)$ has a nonzero element in its second component.
Thus the second component of $q_x+q_y+q_z$ is nonzero,
contradicting our assumption that $q_x+q_y+q_z=0$.
\end{example}

\begin{lemma} \label{lemma:stpp}
If a group $H$ has subsets 
$\{X_i,Y_i,Z_i \,:\, 1 \le i \le n \}$ 
satisfying the simultaneous triple product property,
then for every element $h \pi$ in $H \wr \Sym_n$
there is at most one way to represent $h \pi$ as a
quotient $(x \sigma)^{-1} y \tau$ such that
$x \in \prod_{i=1}^n X_i, \,
y \in \prod_{i=1}^n Y_i, \,
\sigma,\tau \in \Sym_n.$
\end{lemma}
\begin{proof}
Let $X \eqbd \prod_{i=1}^n X_i$, $Y \eqbd \prod_{i=1}^n Y_i.$
Suppose that 
\begin{equation} \label{eqn:lem-stpp-1}
(x \sigma )^{-1} y  \tau = 
(x' \sigma')^{-1} y' \tau' 
\end{equation} 
and that
$x,x' \in X, \, y,y' \in Y, \, 
\sigma,\sigma',\tau,\tau' \in \Sym_n.$
We have
\[
(x \sigma)^{-1} y \tau =
\sigma^{-1} x^{-1} y \tau =
\left (\sigma^{-1} \cdot (x^{-1} y) \right) 
\sigma^{-1} \tau,
\]
and the right side of (\ref{eqn:lem-stpp-1})
may be expressed by a similar formula.
Equating terms on both sides, we find that:
\begin{eqnarray}
\label{eqn:lem-stpp-2}
x_{\sigma(i)}^{-1} y_{\sigma(i)} & = &
(x'_{\sigma'(i)})^{-1} y'_{\sigma'(i)} 
\quad \mbox{for } 1 \le i \le n \\
\label{eqn:lem-stpp-3}
\sigma^{-1} \tau & = & 
(\sigma')^{-1} \tau'.
\end{eqnarray}
Let $j = \sigma(i), k = \sigma'(i).$
We may rewrite
(\ref{eqn:lem-stpp-2}) as
\begin{equation} \label{eqn:lem-stpp-4}
x'_k (x_j)^{-1} y_j (y'_k)^{-1} = 1.
\end{equation}
The left side of (\ref{eqn:lem-stpp-4})
has the form $q_x q_y q_z$ where 
$q_x \in Q(X_k, X_j), \,
q_y \in Q(Y_j, Y_k), \,
q_z = 1 \in Q(Z_k, Z_k).$
By the simultaneous triple product property,
we may conclude that $j=k, \, x=x', \, y=y'.$
Recalling that $j=\sigma(i), \, k=\sigma'(i)$,
we have $\sigma(i)=\sigma'(i)$, and as $i$ was
an arbitrary element of $\{1,2,\ldots,n\}$ we
conclude that $\sigma=\sigma'$.  Combining this
with (\ref{eqn:lem-stpp-3}) implies that
$\tau = \tau'$.  Thus $x=x', \, y=y', \,
\sigma=\sigma', \, \tau=\tau'$, as desired.
\end{proof}

Finally, we must recall some basic facts about
the discrete Fourier transform of an abelian
group.  If $H$ is an abelian group, we let $\widehat{H}$
denote the set of all homomorphisms from $H$ to $S^1$, the
multiplicative group of complex numbers with unit modulus.
Elements of $\widehat{H}$ are called \emph{characters}
and will be denoted in this paper by the letter $\chi$.
The sets $H,\widehat{H}$ have the same cardinality.  
When $H_1, H_2$ are two abelian groups, there is
a canonical bijection between the sets
$\widehat{H_1} \times \widehat{H_2}$
and $(H_1 \times H_2)^{\wedge}$; this bijection
maps an ordered pair $(\chi_1,\chi_2)$ to the
character $\chi$ given by the formula 
$\chi(h_1,h_2) = \chi_1(h_1) \chi_2(h_2).$
Just as the symmetric group $\Sym_n$ acts on $H^n$
via the formula $\sigma \cdot (h_1,h_2,\ldots,h_n)
= (h_{\sigma^{-1}(1)},h_{\sigma^{-1}(2)},\ldots,
h_{\sigma^{-1}(n)}),$ there is a left action of
$\Sym_n$ on the set $\widehat{H}^n$ defined by the
formula $\sigma \cdot (\chi_1,\chi_2,\ldots,\chi_n)
= (\chi_{\sigma^{-1}(1)},\chi_{\sigma^{-1}(2)},\ldots,
\chi_{\sigma^{-1}(n)}).$  In the following section
we will use the notation $\Xi(H^n)$ to denote a
subset of $\widehat{H}^n$ containing exactly one
representative of each orbit of the $\Sym_n$ action
on $\widehat{H}^n$.  An orbit of this action
is uniquely determined by a multiset consisting of 
$n$ characters of $H$, so the cardinality of 
$\Xi(H^n)$ is equal to the number of such multisets,
i.e. $\binom{|H|+N-1}{N}.$

\begin{example}[running example, part 3]
\label{ex_part3}
A character $\chi$ of the group $H = (\Z/16\Z)^3$
is uniquely determined by a triple $(a_1,a_2,a_3)$
of integers modulo $16$.  For an element $h=(b_1,b_2,b_3)
\in H$, we have
$$ \chi(h) = e^{2 \pi i (a_1 b_1 + a_2 b_2 + a_3 b_3) / 16}. $$
A character of the group $H^2$ is a pair of ordered triples
which may be represented as the rows of a matrix
$$ \left( \begin{array}{lll} 
a_{11} & a_{12} & a_{13} \\ a_{21} & a_{22} & a_{23}
\end{array} \right)
$$
as before.  The group $\Sym_2 = \{\pm 1\}$ acts on
$\widehat{H}^2$ by exchanging the two rows of such a
matrix.  An orbit of this action is either:
\begin{itemize}
\item  two distinct matrices, each obtained from the other by
swapping the top and bottom rows; or
\item  one matrix whose top and bottom rows are identical.
\end{itemize}
There are $4096$ rows that can be formed from an ordered
triple of integers modulo $16$, so there are $\binom{4096}{2}$
orbits of the first type and $4096$ orbits of the second type.
Thus the set $\Xi(H^2)$ has cardinality 
$$\binom{4096}{2} + 4096 = 8,\!390,\!656.$$  
\end{example}

\subsection{Abelian STP algorithms}
\label{sec_abelian_stp_alg}

This section is based on the material from~\cite{CKSU-companion}.

\begin{definition}[abelian STP family]
\label{def_stpp_alg}
An \emph{abelian STP family} with growth parameters
$(\alpha,\beta)$ is a collection
of ordered triples $(H_N,\Upsilon_N, k_N)$, defined
for all $N > 0$, satisfying
\begin{enumerate}
\item $H_N$ is an abelian group.
\item 
$\Upsilon_N = {(X_i,Y_i,Z_i} \,:\, i = 1,2,\ldots,N\}$
is a collection of $N$ ordered triples of subsets
of $H_N$ satisfying the simultaneous triple
product property.
\item
$|H_N| = N^{\alpha + o(1)}$.
\item
$k_N = \prod_{i=1}^N |X_i| =
\prod_{i=1}^N |Y_i| =
\prod_{i=1}^N |Z_i| =
N^{\beta N + o(N)}$.
\end{enumerate}
\end{definition}

\begin{remark}
If $\{(H_N,\Upsilon_N,k_N)\}$ is an
abelian STP family, then 
Lemma~\ref{lemma:stpp} ensures
that there is a one-to-one mapping
\[
\left( \prod_{i=1}^N X_i \right)
\times \left( \prod_{i=1}^N Y_i \right) 
\times (\Sym_N)^2 \rightarrow H_N \wr \Sym_N
\]
given by $(x,y,\sigma,\tau) \mapsto (x \sigma)^{-1} y \tau$.
The fact that the mapping is one-to-one implies the first 
line in the following series of inequalities.
\begin{eqnarray*}
|H_N|^N N! & \ge & (k_N N!)^2 \\
N^{\alpha N + o(N)} N^{N+o(N)} & \ge & 
N^{2 \beta N + o(N)} N^{2N + o(N)} \\
\alpha + 1 & \ge & 2 \beta+2 \\
\frac{\alpha-1}{\beta} & \ge & \frac{\alpha+1}{\beta+1} \;\; \ge \;\; 2.
\end{eqnarray*}
\label{remark_alpha_beta}
\end{remark}

\begin{example}[running example, part 4]
\label{ex_part4}
Example~\ref{ex_part3} contained an example of 
a group $H$ with $4096$ elements which contained
two triples of subsets, $(\overline{X}_0,\overline{Y}_0,
\overline{Z}_0)$ and $(\overline{X}_1,\overline{Y}_1,
\overline{Z}_1)$, satisfying the simultaneous triple
product property.  Each of the sets $\overline{X}_i,
\overline{Y}_i,\overline{Z}_i \, (i=0,1)$ has  $15$ 
elements.

We will now show how to extend this
example to an abelian STP family.  For 
$N \ge 1$ let $\ell = \lceil \log_2(N) \rceil$
and let $H_N = H^\ell$.  For $1 \le i \le N$
let $i_1, i_2, \ldots, i_\ell$ denote the
binary digits of the number $i-1$ (padded
with initial $0$'s so that it has exactly
$\ell$ digits) and let
$$ 
X_i \eqbd \prod_{m=1}^\ell \bar{X}_{i_m}, \qquad
Y_i \eqbd \prod_{m=1}^\ell \bar{Y}_{i_m}, \qquad
Z_i \eqbd \prod_{m=1}^\ell \bar{Z}_{i_m}.
$$
The triples $(X_i,Y_i,Z_i)$ satisfy the
simultaneous triple product property.
Indeed, if $i,j,k \in \{1,2,\ldots,N\}$
and $q_x \in Q(X_i,X_j), q_y \in Q(Y_j,Y_k),
q_z \in Q(Z_k,Z_i), q_x + q_y + q_z = 0,$
then for $m=1,2,\ldots,\ell$ it must be
the case that the $m$-th components of
$q_x,q_y,q_z$ satisfy $(q_x)_m + (q_y)_m +
(q_z)_m = 0$.  Using this equation and applying
the fact that $(\overline{X}_0,\overline{Y}_0,
\overline{Z}_0)$ and $(\overline{X}_1,\overline{Y}_1,
\overline{Z}_1)$ satisfy the simultaneous triple
product property, we find that $i_m=j_m=k_m$
and that $(q_x)_m=(q_y)_m=(q_z)_m=0$.  Since this
holds for $m=1,2,\ldots,\ell$, it follows that
$i=j=k$ and $q_x=q_y=q_z=0$ as claimed.

Finally, we may work out the growth parameters
of this abelian STP family.  We have 
$$
|H_N| = |H|^{\ell} = (16^3)^{1 + \lfloor \log_2(N) \rfloor}
= N^{3 \log_2(16) + O(1/\log N)},
$$ 
hence $\alpha = 3 \log_2(16) = 12.$  Also,
$$
k_N = \prod_{i=1}^N |X_i| = 
\prod_{i=1}^N \prod_{m=1}^{\ell} |\bar{X}_{i_m}| =
15^{N \ell} = 15^{N \log_2(N) + O(N)} = N^{N \log_2(15) + O(N/\log N)},
$$
hence $\beta = \log_2(15).$
\end{example}
\bigskip
Given an abelian STP family, one may define a
recursive matrix multiplication which we now
describe.  Given a pair of $n$-by-$n$ matrices
$A,B$, we first find the minimum $N$ such that 
$k_N \cdot N! \ge n$, and we denote the group
$H_N$ by $H$.  If $N! \ge n$, then
we multiply the matrices using an arbitrary
algorithm.  (This is the base of the recursion.)
Otherwise we will reduce the problem
of computing the matrix product $AB$
to $\binom{|H|+N-1}{N}$ instances of 
$N! \times N!$ matrix multiplication,
using a reduction based on the discrete
Fourier transform of the abelian group $H^N$.
We next describe this reduction.

Padding the matrices with additional rows and
columns of $0$'s if necessary, we may assume
without
loss of generality that $k_N \cdot N! = n$.
Define subsets $X,Y,Z \subseteq H \wr \Sym_N$
as follows:
$$
X  \eqbd  \left( \prod_{i=1}^N X_i \right) \times \Sym_N, \qquad
Y  \eqbd  \left( \prod_{i=1}^N Y_i \right) \times \Sym_N, \qquad
Z  \eqbd  \left( \prod_{i=1}^N Z_i \right) \times \Sym_N. 
$$
These subsets satisfy the triple product 
property~\cite{CKSU}.
Note that $|X|=|Y|=|Z|=n.$  We will treat the rows
and columns of $A$ as being indexed by the sets 
$X,Y,$ respectively.  We will treat the rows and 
columns of $B$ as being indexed by the sets $Y,Z,$
respectively.

The algorithm makes use of two auxiliary vector
spaces $\C[H \wr \Sym_N], \C[\widehat{H}^N \rtimes \Sym_N]$, 
each of dimensionality $|H|^N N!$ and each having a
basis which we now designate.  The basis for 
$\C[H \wr \Sym_N]$ is denoted 
by $\{\mathbf{e}_g \,:\, g \in H \wr \Sym_N\}$,
and the basis for
$\C[\widehat{H}^N \rtimes \Sym_N]$ is denoted
by $\{\mathbf{e}_{\chi,\sigma} \,:\, \chi \in \widehat{H}^N,
\sigma \in \Sym_N\}.$

The abelian STP algorithm performs the following
series of steps.  We have labeled the steps according to whether
they perform arithmetic or not. (For example, a step which  permutes 
the components of a vector does not perform arithmetic.)
\begin{enumerate}
\item   \label{step:embed}
\textbf{Embedding}
\textsc{(no arithmetic):\hspace{5mm}}
Compute the following pair of vectors in
$\C[H \wr \Sym_N]$.
\begin{eqnarray*}
a & \eqbd & \sum_{x \in X} \sum_{y \in Y} 
A_{xy} \mathbf{e}_{x^{-1} y} \\
b & \eqbd & \sum_{y \in Y} \sum_{z \in Z} 
B_{yz} \mathbf{e}_{y^{-1} z}. 
\end{eqnarray*}
\item   \label{step:fourier}
\textbf{Fourier transform}
\textsc{(arithmetic):\hspace{5mm}}
Compute the following pair of vectors in
$\C[\widehat{H}^N \rtimes \Sym_N]$.
\begin{eqnarray*}
\hat{a} & \eqbd &
\sum_{\chi \in \widehat{H}^N}
\sum_{\sigma \in \Sym_N} \left(
\sum_{h \in H^N} \chi(h) a_{\sigma h} 
\right) \mathbf{e}_{\chi, \sigma}. \\
\hat{b} & \eqbd &
\sum_{\chi \in \widehat{H}^N}
\sum_{\sigma \in \Sym_N} \left(
\sum_{h \in H^N} \chi(h) b_{\sigma h} 
\right) \mathbf{e}_{\chi, \sigma}. 
\end{eqnarray*}
\item   \label{step:assemble}
\textbf{Assemble matrices}
\textsc{(no arithmetic):\hspace{5mm}}
For every $\chi \in \Xi(H^N)$,
compute the following pair of matrices
$A^\chi, 
B^\chi 
$, whose rows and columns
are indexed by elements of $\Sym_N$.
\begin{eqnarray*}
A^\chi_{\rho \sigma} & \eqbd & \hat{a}_{\rho \cdot \chi, 
\sigma \rho^{-1}} \\
B^\chi_{\sigma \tau} & \eqbd & \hat{b}_{\sigma \cdot \chi,
\tau \sigma^{-1}} \\
\end{eqnarray*}
\item  \label{step:multiply}
\textbf{Multiply matrices}
\textsc{(arithmetic):\hspace{5mm}}
For every $\chi \in \Xi(H^N)$,
compute the matrix product $C^\chi \eqbd A^\chi B^\chi$
by recursively applying the abelian STP algorithm.
\item  \label{step:disassemble}
\textbf{Disassemble matrices}
\textsc{(no arithmetic):\hspace{5mm}}
Compute a vector $\hat{c} \eqbd
 \sum_{\chi,\sigma} \hat{c}_{\chi,\sigma} \mathbf{e}_{\chi,\sigma}
\in \C[\widehat{H}^N \rtimes \Sym_N]$
whose components $\hat{c}_{\chi,\sigma}$ are defined as follows.
Given $\chi,\sigma,$ let $\chi_0 \in \Xi(H^N)$ and $\tau \in
\Sym_N$ be such that $\chi = \tau \cdot \chi_0.$  Let
$$
\hat{c}_{\chi,\sigma} \eqbd C^{\chi_0}_{\tau,\sigma \tau}.
$$
\item  \label{step:inverseFT}
\textbf{Inverse Fourier transform}
\textsc{(arithmetic):\hspace{5mm}}
Compute the following vector $c \in \C[H \wr \Sym_N]$.
$$
c \eqbd \sum_{h \in H^N} \sum_{\sigma \in \Sym_N} \left(
\frac{1}{|H|^N} \sum_{\chi \in \widehat{H}^N} \chi(-h)
\hat{c}_{\chi,\sigma} \right) \mathbf{e}_{\sigma h}.
$$
\item  \label{step:output}
\textbf{Output}
\textsc{(no arithmetic):\hspace{5mm}}
Output the matrix $C = (C_{xz})$ whose entries are 
given by the formula $$C_{xz} \eqbd c_{x^{-1}z}.$$
\end{enumerate}
See~\cite{CKSU-companion} for a proof of the algorithm's
correctness.

\begin{example}[running example, part 5]
\label{ex_part5}

In our example with $H=(\Z/16\Z)^3$ and $N=2$, we have
$k_N N! = (15^2)(2!) = 450$, so the seven
steps outlined above constitute a reduction from 
$450$-by-$450$ matrix multiplication to a large
number of $2$-by-$2$ matrix multiplication problems,
i.e. $|\Xi(H^2)|$ of them.  We will elaborate on the
details of this reduction in the following paragraph. 
Recall from Example~\ref{ex_part3}
that $|\Xi(H^2)| = 8,\!390,\!656$.  By 
comparison, the naive reduction from $450$-by-$450$ to
$2$-by-$2$ matrix multiplication --- by partitioning
each matrix into $(225)^2$ square blocks of size $2$-by-$2$ ---
requires the algorithm to compute $(225)^3 = 11,\!390,\!625$ 
smaller matrix products.  If we use this more efficient 
$450$-by-$450$ matrix multiplication algorithm as the 
recursive step in a stationary partition algorithm as
in Section~\ref{sec_stat}, the running time would be
$O(n^{2.95}).$  Instead, if we use the $N=2, H=(\Z/16\Z)^3$
construction as the basis of an abelian STP family as in 
Example~\ref{ex_part4}, we may apply the abelian STP 
algorithm which uses a more sophisticated recursion as
the size of the matrices grows to infinity.  For example,
when $N=2^\ell$, we have $n = k_N N! = 15^{N\ell} (2^\ell)!$.
The first step of the stationary partition algorithm would
reduce an $n$-by-$n$ matrix multiplication problem to a
set of $(n/450)$-by-$(n/450)$ matrix multiplication problems.
By comparison, the first three steps of the abelian STP algorithm 
reduce $n$-by-$n$ matrix multiplication to a set of $(N!)$-by-$(N!)$
matrix multiplications.  As $N! = O(n^{0.21})$ in this 
example, we see that the abelian STP algorithm achieves 
a much more significant reduction in the size of the matrices
at the top level of recursion.  For the abelian STP algorithm
in our running example, it can be shown that the
running time is $O(n^{2.81})$.

We will now go into greater detail in explaining the
abelian STP algorithm in the case $N=2, H=(\Z/16\Z)^3$
given in our running example.  In this case,
$H \wr \Sym_2$ is the wreath product group described
in Example~\ref{ex_part1}; its elements are represented
by ordered pairs $(M,q)$ where $M$ is a $2$-by-$3$
matrix of integers modulo $16$ and $q = \pm 1.$  The
sets $X,Y,Z \subseteq H \wr \Sym_2$ can be represented
as follows:
\begin{eqnarray*}
X & = & \left\{ \left( \begin{array}{rrr}
\ne 0 & 0 & 0 \\ 0 & \ne 0 & 0 \end{array} \right) \right\} 
\;\; \times \;\; \{\pm 1\} \\[2mm]
Y & = & \left\{ \left( \begin{array}{rrr}
0 & \ne 0 & 0 \\ 0 & 0 & \ne 0 \end{array} \right) \right\} 
\;\; \times \;\; \{\pm 1\} \\[2mm]
Z & = & \left\{ \left( \begin{array}{rrr}
0 & 0 & \ne 0 \\ \ne 0 & 0 & 0 \end{array} \right) \right\} 
\;\; \times \;\; \{\pm 1\}.
\end{eqnarray*}
By this, we mean that an element of $X$ is an ordered pair
$(M,q)$ where  $M$ contains nonzero numbers in the upper
right and lower middle entries, and zero in every other
entry, and $q$ is in $\{\pm 1\}$.  The interpretation of
the expressions for $Y$ and $Z$ is analogous.

The first three steps of the algorithm perform 
preprocessing on the matrix $A$ to arrange some
linear combinations of its entries into a set
of $2$-by-$2$ matrices, one for each element of
$\Xi(H^2).$  Likewise, they preprocess the matrix
$B$ to arrange linear combinations of its 
entries into a set of $2$-by-$2$ matrices.
We will describe the preprocessing of $A$;
the preprocessing of $B$ is entirely analogous,
but uses the subsets $Y,Z \subseteq H \wr \Sym_2$
in place of $X,Y$.
The group $H \wr \Sym_2$ can be partitioned into two
subsets of size $|H|^2=16^6$, namely the elements 
$(M,q)$ whose second component is $+1$ and those whose
second component is $-1$.  (We will call these the
\emph{positive} and \emph{negative} subsets.)  
The first step in the preprocessing of $A$ inserts 
its entries into two $6$-dimensional arrays
$a_+, \, a_-$
of size $16^6$, which we call the \emph{positive} and 
\emph{negative} arrays, corresponding to the positive
and negative subsets of $H \wr \Sym_2$.
For example, the matrix $A$ contains an entry in
row $x = \left( \left( \begin{array}{lll} 
9 & 0 & 0 \\ 0 & 5 & 0 \end{array} \right), \;
-1 \right)$ and column 
$y = \left( \left( \begin{array}{lll}
0 & 11 & 0 \\ 0 & 0 & 4 \end{array} \right), \;
1 \right)$, because $x \in X$ and $y \in Y$.
In $H \wr \Sym_2$ we may compute that 
$$
x^{-1} y = \left( \left( \begin{array}{lll}
0 & 11 & 0 \\ 0 & 0 & 4 \end{array} \right) \,-\,
\left( \begin{array}{lll}
0 & 5 & 0 \\ 9 & 0 & 0 \end{array} \right), \;
-1 \right) = \left( \left( \begin{array}{lll}
0 & 6 & 0 \\ 7 & 0 & 4 \end{array} \right), \;
-1 \right).
$$
This tells us that the entry $A_{xy}$ in row $x$,
column $y$ of $A$ should be inserted into $a_-$
(because the second component of
$x^{-1} y$ is -1) at the location whose index
in the 6-dimensional array is the 6-tuple
$(0,6,0,7,0,4).$  The locations of the other entries
of $A$ are determined similarly.  At the end of this
step, some of the entries of the positive and negative
arrays will not have been filled with an entry of $A$;
the algorithm writes $0$ in these entries of the 
positive and negative arrays.

The second step in the preprocessing of $A$ 
performs a 6-dimensional discrete Fourier transform on
the positive and negative arrays.  That is, we compute
the array $\hat{a}_{+}$ whose entries are given
by the formula:
\[
\hat{a}_{+}(i_1,i_2,i_3,i_4,i_5,i_6) =
\sum_{j_1 = 0}^{15}
\sum_{j_2 = 0}^{15}
\sum_{j_3 = 0}^{15}
\sum_{j_4 = 0}^{15}
\sum_{j_5 = 0}^{15}
\sum_{j_6 = 0}^{15}
\exp \left( \frac{2 \pi i}{16} \sum_{k=1}^t i_k j_k \right)
a(j_1,j_2,j_3,j_4,j_5,j_6).
\]
This may be computed using the fast Fourier transform.
An array $\hat{a}_{-}$ is defined similarly, using the
entries of the negative array instead of the positive
array.

The third step in the preprocessing of $A$ forms
a $2$-by-$2$ matrix $A^\chi$ for each element $\chi \in \Xi(H^2)$.
The formula is given in step~\ref{step:assemble} above.
Recall that an element of $\Xi(H^2)$ can be represented by
a $2$-by-$3$ matrix of integers modulo $16$, subject to
the condition that if two such matrices differ only by
swapping the top and bottom rows, then exactly one of
them belongs to $\Xi(H^2)$.  For notational convenience,
we will write the entries of a matrix 
$\left( \begin{array}{lll} i_1 & i_2 & i_3 \\ i_4 & i_5 & i_6
\end{array} \right)$ as a 6-tuple $(i_1,i_2,i_3,i_4,i_5,i_6).$
If $\chi=(i_1,i_2,i_3,i_4,i_5,i_6)$ then $A^\chi$ is the matrix
\[
\left( \begin{array}{ll}
\hat{a}_+(i_1,i_2,i_3,i_4,i_5,i_6) &
\hat{a}_-(i_1,i_2,i_3,i_4,i_5,i_6) \\
\hat{a}_-(i_4,i_5,i_6,i_1,i_2,i_3) &
\hat{a}_+(i_4,i_5,i_6,i_1,i_2,i_3)
\end{array} \right).
\]
Note that if $i_1=i_4, i_2=i_5, i_3=i_6$ then this matrix
contains only two distinct numbers, each repeated twice.
The preprocessing of $B$ is performed similarly, resulting
in matrices $B^\chi$ for each $\chi \in \Xi(H^2)$.  The
algorithm then computes each matrix product 
$C^\chi = A^\chi B^\chi.$  

Finally, there is a three-step
postprocessing phase which reconstructs the entries of the
matrix product $C = AB$ by taking linear combinations of the
entries of the matrices $C^\chi$.  The first step is to
arrange the entries of the matrices $C^\chi$ into a pair
of arrays $\hat{c}_+, \, \hat{c}_-$ by reversing the mapping
which was used to assemble the entries of $\hat{a}_+, \,
\hat{a}_-$ into the matrices $A^\chi$.  Thus, for a 
6-tuple $\chi \eqbd (i_1,i_2,i_3,i_4,i_5,i_6)$, if $\chi \in \Xi(H^2)$
then $\hat{c}_+(i), \hat{c}_-(i)$ are the entries of the
first row of $C^\chi$ and if $i \not\in \Xi(H^2)$ then
$\hat{c}_-(i), \hat{c}_+(i)$ are the entries of the
second row of $C^{\chi'}$ where $\chi' \eqbd (i_4,i_5,i_6,i_1,i_2,i_3).$
Having constructed the arrays $\hat{c}_+, \, \hat{c}_-$, we
perform an inverse Fourier transform to obtain arrays
$c_+, c_-$.  Finally, to determine the entry $C_{xz}$
of the product matrix $C = AB$, we compute the element
$x^{-1} z$ in the wreath product $H \wr \Sym_2$, 
select the array $c_+$ or $c_-$ according to whether
the second component of $x^{-1} z$ is $+1$ or $-1$, 
and look up the entry in this array whose index is
the 6-tuple given by the entries of the matrix which
forms the first component of $x^{-1} z$.
\end{example}

\subsection{Analysis of abelian STP algorithms} \label{sec_grtheory_analysis}

Now we are in a position to analyze abelian STP algorithms. We could have done
that using Theorem~\ref{thm_crude}. However, we choose to further refine our error 
analysis 
to obtain sharper norm inequalities for a specific matrix norm and hence better error bounds.

\begin{theorem} \label{thm_error_analysis}
If $\{(H_N,\Upsilon_N,k_N ) \}$ is an abelian
STP family with growth parameters $(\alpha,\beta)$,
then the corresponding abelian STP algorithm
is stable. It satisfies the error bound~(\ref{gen_bound}), 
with the Frobenius norm and the function $\mu$ of order 
$$  \mu(n)=n^{\frac{\alpha+2}{2\beta} \,+\, o(1)}. $$
\end{theorem}
\begin{proof}
We seek to establish that the matrix $C_{comp}$ computed by the algorithm
differs from the actual matrix product $C$ by at most
$\mu(n) \ve \Fr{A} \Fr{B}+O(\ve^2)$ in Frobenius norm, i.e.
\[
\Fr{C_{comp} - C} \le \mu(n) \ve \Fr{A} \Fr{B} +O(\ve^2).
\]

Throughout this proof we will adopt the
convention that the Fourier transform of
an abelian group $H$ is represented by a
matrix $F$ satisfying $F F^{\mathsf{T}} 
= |H|I$, rather than a unitary matrix.
This is consistent with the interpretation
that the Fourier transform takes an
element $x \in \C[H]$, represented as
a linear combination of basis elements
$h \in H$, and returns the coefficients
$a_\chi$ in the unique representation
of $x$ as a linear combination
$\sum_\chi a_\chi w_\chi,$ where the elements
$w_\chi$ are idempotent in $\C[H].$  When the
Fourier transform is instead normalized
so that it is represented by a unitary
matrix, we will refer to this linear
transformation as the ``unitary Fourier
transform.''
Let $f(n)$ be the $L_2$ error bound 
satisfied by the unitary Fourier transform
and its inverse, i.e. if $H$ is an abelian group
with $n$ elements, $x$ is a vector in
$\C[H]$, and $\hat{x}, \check{x}$ are
its unitary Fourier transform and inverse Fourier
transform, then $\|\hat{x}_{comp}-\hat{x}\|_2$
and $\|\check{x}_{comp}-\check{x}\|_2$
are both bounded by $f(n) \ve \|x\|_2 + O(\ve^2)$.

An abelian STP algorithm satisfies a Frobenius-norm
error bound of the form $\mu(n)$, where the function
$\mu(n)$ satisfies a recursion which is governed by
the recursive structure of the algorithm itself.  
Specifically, the algorithm break down into 
a series of seven steps specified in 
Section~\ref{sec_abelian_stp_alg}.
Observe that arithmetic operations are performed
only in the even-numbered steps.  
The odd-numbered steps consist only of 
rearranging (and possibly repeating)
the components of a vector
to form the entries of a set of matrices
and vice-versa.  (To see that
no arithmetic is performed in Step~\ref{step:embed},
use Lemma~\ref{lemma:stpp} which implies that 
each component of the vectors $a,b$ is a sum of either
zero or one entry of one of the matrices $A,B$.)

Step~\ref{step:embed} replaces the matrix
$A$ with a vector whose $2$-norm is
equal to the Frobenius norm $\Fr{A}$,
and similarly for $B$.  Step~\ref{step:fourier}
performs $N!$ copies of the discrete
Fourier transform of the group $H^N$.
We have
\begin{equation}
\|\hat{a}_{comp}-\hat{a}\|_2^2 \le
|H|^N f(|H|^N)^2 \ve^2 \Fr{A}^2 \;+\; O(\ve^3)
\end{equation}
and similarly for $B$.  (The extra
factor of $|H|^N$ on the right side
comes from the fact that we're using 
a Fourier transform which is a unitary
matrix multiplied by the scalar 
$|H|^{N/2}.$)

Step~\ref{step:assemble} doesn't perform any
arithmetic, but it repeats each component
of $\hat{a}$ (resp. $\hat{b}$) possibly
$N!$ times in assembling a set of
matrices $\{A^\chi\}$ (resp.$\{B^\chi\}$).
Let $\Err_A^\chi \eqbd A^\chi_{comp} - A^\chi$.
We have
\begin{eqnarray} 
\label{eqn:Achi}
\sum_{\chi \in \Xi} \Fr{A^\chi}^2 & \le &
N! |H|^N \Fr{A}^2 \\
\label{eqn:errA}
\sum_{\chi \in \Xi}
\Fr{\Err_A^\chi}^2 & \le &
 N! |H|^N f(|H|^N)^2 \ve^2 \Fr{A}^2 \;+\; O(\ve^3).
\end{eqnarray}
The matrices $\Err_B^\chi$ are defined
similarly, and they satisfy a similar
bound on the sum of their squared Frobenius
norms.

Step~\ref{step:multiply} multiplies each pair
$A^\chi_{comp}, B^\chi_{comp}$ to obtain a
matrix $C^\chi_{comp}$.  The error matrix
$$
\Err_C^\chi \eqbd C^\chi_{comp} - C^\chi
$$
can be expressed as a sum of four terms, as
follows:
\begin{eqnarray*}
\Err_C^\chi & = & (C^\chi_{comp} -
A^\chi_{comp} B^\chi_{comp}) +
(A^\chi_{comp} B^\chi_{comp} - C_\chi) \\
& = &
(C^\chi_{comp} - A^\chi_{comp} B^\chi_{comp})
+ [(A^\chi + \Err_A^\chi)(B^\chi + \Err_B^\chi)
- A^\chi B^\chi] \\
& = &
(C^\chi_{comp} - A^\chi_{comp} B^\chi_{comp})
+ \Err_A^\chi B^\chi + A^\chi \Err_B^\chi
+ \Err_A^\chi \Err_B^\chi.
\end{eqnarray*}
The fourth term is of order $O(\ve^2)$ and
may be ignored.  The remaining terms may be
dealt with as follows.  First, by the inductive
hypothesis:
\begin{eqnarray*}
\Fr{C^\chi_{comp} - A^\chi_{comp} B^\chi_{comp}}
& \le & \mu(N!) \ve \Fr{A^\chi_{comp}} \Fr{B^\chi_{comp}} \;+\; O(\ve^2)\\
& = & \mu(N!) \ve \Fr{A^\chi} \Fr{B^\chi} \;+\; O(\ve^2).
\end{eqnarray*}
Next,
\begin{eqnarray*}
\Fr{\Err_A^\chi B^\chi} & \le & \Fr{\Err_A^\chi} \Fr{B^\chi} \\
\Fr{A^\chi \Err_B^\chi} & \le & \Fr{A^\chi} \Fr{\Err_B^\chi}.
\end{eqnarray*}
Summing all of these bounds over $\chi \in \Xi$, we obtain
\begin{eqnarray*}
\sum_\chi \Fr{\Err^C_\chi} & \le &
\mu(N!) \ve \sum_\chi \Fr{A^\chi} \Fr{B^\chi} +
\sum_\chi \Fr{\Err_A^\chi} \Fr{B^\chi} +
\sum_\chi \Fr{A^\chi} \Fr{\Err_B^\chi} \;+\; O(\ve^2) \\
& \le &
\mu(N!) \ve 
\left(\sum_\chi \Fr{A^\chi}^2 \right)^{1/2}
\left(\sum_\chi \Fr{B^\chi}^2 \right)^{1/2} + \\
& & \quad \left(\sum_\chi \Fr{\Err_A^\chi}^2 \right)^{1/2}
\left(\sum_\chi \Fr{B^\chi}^2 \right)^{1/2} + \\
& & \quad
\left(\sum_\chi \Fr{A^\chi}^2 \right)^{1/2}
\left(\sum_\chi \Fr{\Err_B^\chi}^2 \right)^{1/2} \;+\;
O(\ve^2) \\
& \le &
\mu(N!) N! |H|^N \ve \Fr{A} \Fr{B} +
2 f(|H|^N) N! |H|^N  \ve \Fr{A} \Fr{B} \;+\; O(\ve^2) \\
\left(\sum_\chi \Fr{\Err^C_\chi}^2\right)^{1/2} & \le &
\left[ 2 f(|H|^N) + \mu(N!) \right]
N! |H|^N \ve \Fr{A} \Fr{B} \;+\; O(\ve^2) .
\end{eqnarray*}
The second line was derived from the first
by applying Cauchy-Schwarz three times.  The
third line was derived using (\ref{eqn:Achi})
and (\ref{eqn:errA}).
The final line was derived using the 
inequality $\|x\|_2 \le \|x\|_1$,
applied to the vector whose components
are $(\Fr{\Err^C_\chi})_{\chi \in \Xi}$.

In Step~\ref{step:disassemble}, we form a vector $\hat{c}_{comp}$ in
Fourier space whose components are a subset of
the entries of the matrices $C^\chi_{comp}$.
Our upper bound on 
$\left(\sum_\chi \Fr{\Err^C_\chi}^2\right)^{1/2}$
remains a valid upper bound on
$\|\hat{c}_{comp}-\hat{c}\|_2$.

In Step~\ref{step:inverseFT}, we apply the inverse 
Fourier transform to $\hat{c}_{comp},$ to obtain a
vector $c_{comp}.$  
The inverse Fourier transform performed in this step 
is a unitary Fourier transform 
multiplied by the scalar $|H|^{-N/2},$ so
\begin{eqnarray*}
\|c_{comp}-c\|_2 & \le &
|H|^{-N/2} f(|H|^N) \ve \|\hat{c}_{comp}\|_2 +
|H|^{-N/2} \|\hat{c}_{comp} - \hat{c}\|_2 \;+\; O(\ve^2) \\
& \le &
|H|^{-N/2} f(|H|^N) \ve \|\hat{c}\|_2 +
\left[2 f(|H|^N) + \mu(N!)\right] N! |H|^{N/2}
\ve \Fr{A} \Fr{B} \;+\; O(\ve^2) \\
& \le &
f(|H|^N) \ve \Fr{A} \Fr{B} +
\left[2 f(|H|^N) + \mu(N!) \right] N! |H|^{N/2}
\ve \Fr{A} \Fr{B} \;+\; O(\ve^2).
\end{eqnarray*}
The second line was derived from the first by
observing that $\|\hat{c}_{comp}\|_2 \le \|\hat{c}\|_2+O(\ve)$
and by substituting our earlier bound for
$\|\hat{c}_{comp} - \hat{c}\|_2.$  The third line
was derived by using the bound 
$\|\hat{c}\|_2 \le |H|^{N/2} \Fr{A} \Fr{B},$
which follows from the fact that $\hat{c}$ is the Fourier
transform of the vector $c$, whose $L_2$-norm is
$\Fr{C} \le \Fr{A} \Fr{B}.$  (Recall that the Fourier 
transform increases $L_2$ norms of vectors by a factor 
of $|H|^{N/2}.$)

In Step~\ref{step:output}, no further error is introduced.
Thus, the matrix $C_{comp}-C$ has Frobenius
norm bounded by
\begin{equation} \label{eqn:error-final}
\Fr{C_{comp}-C} \le \left[  
f(|H|^N) + 2 N! |H|^{N/2} f(|H|^N) + N! |H|^{N/2} \mu(N!) \right]
\ve \Fr{A} \Fr{B} \;+\; O(\ve^2).
\end{equation}
Let $m=N!$, and recall that $n=m^{\beta+1+o(1)},
|H|^N = m^{\alpha+o(1)}.$
The error bound (\ref{eqn:error-final}) leads to
the recursion
\[
\mu(m^{\beta+1+o(1)}) \le
m^{1+\alpha/2+o(1)} f(m^{\alpha+o(1)}) + m^{1+\alpha/2+o(1)} \mu(m).
\]
Assuming that the Fourier transform is implemented
using the Cooley-Tukey FFT (see, e.g.,~\cite{MaslenRockmore}), we have $f(n)=O(\log n)$.
Solving the recursion, we find that
\[
\mu(m) = m^{\frac{\alpha+2}{2\beta} \,+\, o(1)}.
\]
\end{proof}

\begin{remark} 
Note that we could apply Theorem~\ref{thm_crude} directly, with the pre-processing map
performing steps~\ref{step:embed} through~\ref{step:assemble} of the algorithm, and the post-processing map performing 
steps~\ref{step:disassemble} through~\ref{step:output}. 
From the discussion in this section we see that the operator norms
of these maps subordinate to the Frobenius norm are bounded as
$$ \| \pre_n \|_{op} \leq (N!)^{1/2} |H|^{N/2}, \qquad \| \post_n \|_{op} \leq |H|^{-N/2}, $$
while the error functions $f_{pre}$ and $f_{post}$ are bounded by $$ f_{pre}(n) \leq  
(N!)^{1/2} |H|^{N/2}  f(|H|^N), \qquad f_{post}(n) \leq |H|^{-N/2} f(|H|^N).$$
Finally, the number of blocks $t$ is ${|H|+N-1\choose N} \approx |H|^N / N!$.   
This leads to a bound 
\begin{equation} \label{eqn:error-finalcrude}
\Fr{C_{comp}-C} \le \left[|H|^{3N/2} \mu(N!)+ 2 |H|^{N} (N!)^{-1/2}  f(|H|^N) + 
 N! |H|^{N/2} f(|H|^N) \right] \ve \Fr{A} \Fr{B}+O(\ve^2),
\end{equation}
which is somewhat weaker than~(\ref{eqn:error-final}). From~(\ref{eqn:error-finalcrude}),
we then obtain the recursion
\[ \mu(m^{\beta+1+o(1)}) \le m^{3\alpha/2+o(1)} \mu(m) +
m^{1+\alpha/2+o(1)} f(m^{\alpha+o(1)}), \]
which gives $$\mu(m) = m^{\frac{3\alpha}{2\beta} +o(1)}.$$
\end{remark}

\begin{remark}
The running time of an abelian STP algorithm can also be
bounded in terms of the growth parameters of the abelian
STP family.  Specifically, the running time is 
$O \left( n^{(\alpha-1)/\beta \, + \, o(1)} \right).$
See~\cite{CKSU-companion} for details. 
Note that the sum of the two exponents, $(\alpha-1)/\beta$ and
$(\alpha+2)/2\beta$, is always bigger than $3$, since
$\alpha \geq 2\beta+1$: 
$$ {\alpha -1\over \beta} + {\alpha +2 \over 2\beta}={3\alpha\over 2\beta}
\geq {6\beta+3\over 2\beta}>3 .$$

\end{remark}

\subsection{Analysis of examples} \label{sec_examples}

The abelian STP algorithm in our running example has
growth parameters $\alpha=12, \beta=\log_2(15)$, hence
its running time is 
$$
O \left( n^{\frac{\alpha-1}{\beta} \,+\, o(1)} \right) = 
O \left( n^{2.82} \right),$$ 
and the error bound is 
$$O \left( n^{\frac{\alpha+2}{2 \beta} \,+\, o(1)} \right) =
O \left( n^{1.80} \right).$$
Note that this error bound (in the Frobenius norm) implies
a bound of $O(n^{2.80})$ in the entrywise maximum norm.
This compares favorably with the error bound for Strassen's
algorithm, which is $O(n^{3.58})$ in the entrywise maximum
norm, while nearly matching the exponent in the running time
of Strassen's algorithm, which runs in time $O\left( n^{\log_2(7)}
\right) = O\left(n^{2.81} \right)$.
Many other examples of abelian STP algorithms are
listed in~\cite{CKSU}.  The algorithms with
running time $O(n^{2.48})$  in Propositions 3.8 and
4.5 of~\cite{CKSU} are both based
on an explicit construction of an abelian STP 
family with growth parameters
$\alpha = 3 \log_4(6), \beta = \log_4(5)$,
hence the error bound for these algorithms
is $O(n^{2.54})$.
The algorithm with running time $O(n^{2.41})$ 
described in Theorems 3.3 and 6.6 of~\cite{CKSU} is
based on an abelian STP family with growth
parameters $\alpha = 3 \log_{6.75}(10),
\beta = \log_{6.75}(8),$
hence the error bound for this algorithm is
$O(n^{2.58})$.

\section{Stability of linear algebra algorithms based on matrix multiplication}\label{sec_addla}

It is natural to ask what other linear algebra operations
  can be done stably and quickly by depending on the
  stability of fast matrix multiplication described here.
  Indeed, ``block'' algorithms relying on matrix multiplication
  are used in practice for many linear algebra operations
  \cite{lapackmanual3,scalapackmanual}, and have been
  shown to be stable assuming only the error bound (\ref{gen_bound})
  \cite{DemmelHigham92}. In a companion paper \cite{new},
 we show that while stable these earlier block algorithms are not
 asymptotically as fast as matrix multiplication.
 However, \cite{new} also shows there are variants of these block
 algorithms for operations like QR decomposition, linear equation
 solving and determinant computation that are both stable and
 as fast as matrix multiplication.

\section{Acknowledgements}

We thank Henry Cohn,
Bal\'azs Szegedy, and Chris Umans for helpful
discussions about this work. We acknowledge both
Alicja Smoktunowicz and Doug Arnold for pointing out~\cite{Miller75}.

\bibliographystyle{plain}
\bibliography{biblio,matrix2}

\end{document}